# A generalisation of the Bernoulli numbers from the discrete to the continuous

Donal F. Connon

dconnon@btopenworld.com

16 May 2010

## Abstract

We generalise the Bernoulli numbers to include the case where the index may be a continuous variable.

## 1. Introduction

The Bernoulli polynomials $B_n(x)$ in the case where $n$ is an integer $n \geq 0$ are defined by the series [2, p.264]

(1.1) $$\frac{te^{tx}}{e^t - 1} = \sum_{n=0}^{\infty} B_n(x) \frac{t^n}{n!} \qquad , (|t| < 2\pi)$$

and the Bernoulli numbers $B_n$ are given by the generating function

(1.2) $$\frac{t}{e^t - 1} = \sum_{n=0}^{\infty} B_n \frac{t^n}{n!} \qquad , (|t| < 2\pi)$$

and it is readily seen from (1.1) that $B_n(0) = B_n$.

The Bernoulli numbers were originally introduced by Johann Faulhaber (1580-1635) in his book Academia Algebrae published in 1631. To his credit, Jacob Bernoulli (1654-1705) did in fact give priority to Faulhaber in his treatise on probability theory, Ars Conjectandi, which was published posthumously in 1713.

Using the Cauchy product [9, p.146] of two infinite series, equation (1.1) may be written as follows

$$\sum_{n=0}^{\infty} B_n(x) \frac{t^n}{n!} = \left(\frac{t}{e^t - 1}\right) e^{tx} = \left(\sum_{n=0}^{\infty} B_n \frac{t^n}{n!}\right)\left(\sum_{n=0}^{\infty} x^n \frac{t^n}{n!}\right)$$

$$= \sum_{n=0}^{\infty} t^n \left(\sum_{k=0}^{n} \frac{B_k x^{n-k}}{k!(n-k)!}\right)$$

and equating coefficients of $t^n$ we have

(1.3) $$B_n(x) = \sum_{k=0}^{n} \binom{n}{k} B_k x^{n-k}$$

This identity could also be obtained in a less rigorous manner by formally differentiating (1.1) with respect to $x$. It is therefore evident from (1.3) that $B_n(x)$ is a polynomial of degree $n$ and that the coefficient of $x^n$ is 1 (i.e., $B_n(x)$ is a monomial because $B_0 = 1$).

We have the identity

$$\frac{te^{t(x+1)}}{e^t - 1} - \frac{te^{tx}}{e^t - 1} = te^{tx}$$

from which we find

$$\sum_{n=0}^{\infty} [B_n(x+1) - B_n(x)] \frac{t^n}{n!} = \sum_{n=0}^{\infty} x^n \frac{t^{n+1}}{n!}$$

and equating the coefficients of $t^n$ gives us for $n \ge 1$

(1.4) $$B_n(x+1) - B_n(x) = nx^{n-1}$$

Letting $x = 0$ in (1.4) results in for $n \ge 2$

(1.5) $$B_n(1) = B_n(0) = B_n$$

and using (1.3) we obtain the recurrence relation for $n \ge 2$

(1.6) $$B_n = \sum_{k=0}^{n} \binom{n}{k} B_k$$

We also have an explicit formula for $B_n$ given by F. Lee Cook in [4]. In what follows, I have applied Lee Cook's approach to the Bernoulli polynomial, rather than to $B_n$ .

Using (1.1) and the Maclaurin expansion we have

$$B_n(x) = \left[ \frac{d^n}{dt^n} \left( \frac{te^{tx}}{e^t - 1} \right) \right]_{t=0}$$

Since $t = \log\left[1 - (1 - e^t)\right]$ and $\log(1-u) = -\sum_{m=1}^{\infty} \frac{u^m}{m}$ for $|u| < 1$, we have

$$t = -\sum_{m=1}^{\infty} \frac{(1-e^t)^m}{m}$$

Therefore, dividing by $(e^t - 1)$ we get

$$\frac{te^{tx}}{e^t - 1} = \sum_{m=1}^{\infty} \frac{e^{tx}(1-e^t)^{m-1}}{m} = \sum_{k=0}^{\infty} \frac{e^{tx}(1-e^t)^k}{k+1} \quad \text{for } |t| < \log 2$$

Hence we have

$$B_n(x) = \sum_{k=0}^{\infty} \frac{1}{k+1} \left[ \frac{d^n}{dt^n} \left(1-e^t\right) e^{tx} \right]_{t=0}$$

Using the Leibniz rule for the derivative of a product we note that the $n$ th derivative vanishes at $t = 0$ for $k \geq n+1$. Therefore, using the binomial theorem we have

$$B_n(x) = \sum_{k=0}^{\infty} \frac{1}{k+1} \sum_{j=0}^{k} (-1)^j \binom{k}{j} \left[ \frac{d^n}{dt^n} (e^{(j+x)t}) \right]_{t=0}$$

Hence we have

$$(1.7) \qquad B_n(x) = \sum_{k=0}^{n} \frac{1}{k+1} \sum_{j=0}^{k} (-1)^j \binom{k}{j} (x+j)^n$$

Since $\sum_{j=0}^{k} (-1)^j \binom{k}{j} (x+j)^n = 0$ for $k > n = 0,1,2,...$ we therefore have the infinite series representation

$$(1.8) \qquad B_n(x) = \sum_{k=0}^{\infty} \frac{1}{k+1} \sum_{j=0}^{k} (-1)^j \binom{k}{j} (x+j)^n$$

A different proof, using the Hurwitz-Lerch zeta function, was recently given by Guillera and Sondow [6].

When $x = 0$ we obtain

$$(1.9) \qquad B_n(0) = B_n = \sum_{k=0}^{n} \frac{1}{k+1} \sum_{j=0}^{k} (-1)^j \binom{k}{j} j^n$$

$$= \sum_{k=0}^{\infty} \frac{1}{k+1} \sum_{j=0}^{k} (-1)^j \binom{k}{j} j^n$$

The structural similarity of (1.9) with the Stirling numbers of the second kind [11, p.58] may immediately be noted

(1.10) $$S(n,k) = \frac{1}{k!} \sum_{\mathrm{j}=0}^{k} (-1)^{k-j} \binom{k}{j} j^n$$

We therefore conclude that

(1.11) $$B_n = \sum_{k=0}^{n} \frac{(-1)^k k!}{k+1} S(n,k)$$

Surprisingly, this simple relationship does not appear in the book "Concrete Mathematics" [5]: instead, the authors report a much more complex identity involving both kinds of the Stirling numbers [5, p.289]. I subsequently discovered that the formula (1.11) is reported in the Wolfram Mathworld website dealing with the Bernoulli numbers. Furthermore, this identity was also proved by Kaneko in 2000 in his paper "The Akiyama-Tanigawa algorithm for Bernoulli numbers" [8], albeit his method was somewhat less direct. Interestingly, in that paper Kaneko reports that the Bernoulli numbers were independently discovered by Takakazu Seki (1642-1708) one year prior to Jacob Bernoulli.

Another proof of (1.11) was also given by Akiyama and Tanigawa in [1]. That paper also mentions that

$$a_k = \sum_{j=1}^{k} s(k,j) b_j \quad \Leftrightarrow \quad b_k = \sum_{j=1}^{k} S(k,j) a_j$$

and hence we obtain

(1.12) $$\sum_{r=1}^{k} s(k,r) B_r = \frac{(-1)^k k!}{k+1}$$

## 2. A generalisation of the Bernoulli numbers from the discrete to the continuous

In a curious paper published in 2000 Musès [10] strongly argued that the discrete Bernoulli numbers should be extended to a Bernoulli function which he denoted by $\beta(s)$ (where $s$ may be real or complex).

Using Euler's formula

(2.1) $$B_{2n} = 2(-1)^{n+1}(2n)!\varsigma(2n)(2\pi)^{-2n}$$

and making the substitution $s = 2n$ and noting that

$$(-1)^{n+1} = \cos\left(2n\pi/2\right)\cos[\pi(1-2n)]$$

Musès contended that the Bernoulli function $\beta(s)$ could be represented by

(2.2) $$\beta(s) = 2\Gamma(s+1)\varsigma(s)(2\pi)^{-s}\cos\left(\frac{\pi s}{2}\right)\cos[\pi(1-s)]$$

so that when $s$ is an integer we have

(2.3) $$\beta(n) = B_n$$

Therefore if $s$ is an even integer we have

$$B_{2n} = \beta(2n) = 2(2n)!\varsigma(2n)(2\pi)^{-2n}\cos\left(n\pi\right)\cos[\pi(1-2n)]$$

$$= 2(-1)^{n+1}(2n)!\varsigma(2n)(2\pi)^{-2n}$$

With $s = 2n+1$ Musès reports that

(2.4) $$B_{2n+1} = \beta(2n+1) = 2(2n+1)!\varsigma(2n+1)(2\pi)^{-2n-1}\cos\left(\frac{[2n+1]\pi}{2}\right)$$

whereupon we immediately see that this corresponds with $B_{2n+1} = 0$ for all $n \ge 1$.

Furthermore, if (2.4) is valid as $n \to 0$ we have

$$B_1 = \beta(1) = \frac{1}{\pi}\lim_{n\to 0}\varsigma(2n+1)\cos\left(\frac{[2n+1]\pi}{2}\right)$$

and from Riemann's functional equation [2, p.259]

(2.5) $$\varsigma(1-s) = 2(2\pi)^{-s}\Gamma(s)\cos(\pi s/2)\varsigma(s)$$

we may determine that

(2.5.1) $$\lim_{s\to 1}\varsigma(s)\cos\left(\frac{s\pi}{2}\right) = -\frac{\pi}{2}$$

Hence we may deduce that $\beta(1) = -\frac{1}{2}$ which corresponds with the well known result that $B_1 = -\frac{1}{2}$. The limit in (2.5.1) may also be ascertained from

$$\lim_{s\to1} \varsigma(s)\cos(s\pi/2) = \lim_{s\to1}(s-1)\varsigma(s)\frac{\cos(s\pi/2)}{s-1} = \lim_{s\to1}(s-1)\varsigma(s)\lim_{s\to1}\frac{\cos(s\pi/2)}{s-1}$$

Musès also contended that for negative integers

(2.6) $$B_{-2n} = 2n\varsigma(2n+1)$$

and

(2.7) $$B_{1-2n} = (1-2n)\varsigma(2n)$$

and it may be noted that this is related to (2.21).

We now consider a function $B_s(1)$ defined by

(2.8) $$B_s(1) = \sum_{m=0}^{\infty}\frac{1}{m+1}\sum_{k=0}^{m}\binom{m}{k}(-1)^k(1+k)^s$$

and link this with the Hasse identity [7] for the Hurwitz zeta function which holds for all $s \in \mathbf{C}$ except $s = 1$

(2.9) $$\varsigma(s,a) = \frac{1}{s-1}\sum_{m=0}^{\infty}\frac{1}{m+1}\sum_{k=0}^{m}\binom{m}{k}\frac{(-1)^k}{(k+a)^{s-1}}$$

We have

$$\lim_{s\to1}\left[\varsigma(s,a) - \frac{1}{s-1}\right] = \lim_{s\to1}\frac{1}{s-1}\left[\sum_{m=0}^{\infty}\frac{1}{m+1}\sum_{k=0}^{m}\binom{m}{k}\frac{(-1)^k}{(k+a)^{s-1}} - 1\right]$$

and using L'Hôpital's rule this becomes

$$= -\lim_{s\to1}\sum_{m=0}^{\infty}\frac{1}{m+1}\sum_{k=0}^{m}\binom{m}{k}\frac{(-1)^k\log(k+a)}{(k+a)^{s-1}}$$

$$= -\sum_{m=0}^{\infty}\frac{1}{m+1}\sum_{k=0}^{m}\binom{m}{k}(-1)^k\log(k+a)$$

$$= -\psi(a)$$

With $a = 1$ in (2.9) we have

(2.9.1) $$\varsigma(s) = \frac{1}{s-1} \sum_{m=0}^{\infty} \frac{1}{m+1} \sum_{k=0}^{m} \binom{m}{k} \frac{(-1)^k}{(k+1)^{s-1}}$$

We see that the substitution $s \to 1-s$ gives us

$$\varsigma(1-s) = -\frac{1}{s} \sum_{m=0}^{\infty} \frac{1}{m+1} \sum_{k=0}^{m} \binom{m}{k} (-1)^k (k+1)^s$$

and using Riemann's functional equation (2.5) this may be written as

$$2(2\pi)^{-s} \Gamma(s) \cos(\pi s/2) \varsigma(s) = -\frac{1}{s} \sum_{m=0}^{\infty} \frac{1}{m+1} \sum_{k=0}^{m} \binom{m}{k} (-1)^k (k+1)^s$$

or equivalently

(2.10) $$2(2\pi)^{-s} \Gamma(s+1) \varsigma(s) \cos(\pi s/2) = -\sum_{m=0}^{\infty} \frac{1}{m+1} \sum_{k=0}^{m} \binom{m}{k} (-1)^k (k+1)^s$$

Comparing this with (2.2) we deduce that

(2.11) $$\beta(s) = -\cos[\pi(1-s)] B_s(1)$$

and we see from (1.8) and (2.8) that $B_s(1)$ is equal to $B_n$ when $s = n$. This confirms that

(2.12) $$\beta(n) = -\cos[\pi(1-n)] B_n$$

so that

$$\beta(2n) = B_{2n}$$

and

$$\beta(2n+1) = -B_{2n+1}$$

Letting $s = -n$ in (2.9) gives us

$$\varsigma(-n, a) = -\frac{1}{n+1} \sum_{m=0}^{\infty} \frac{1}{m+1} \sum_{k=0}^{m} \binom{m}{k} (-1)^k (k+a)^{n+1}$$

and referring to (1.8) we immediately see that [2, p.264]

(2.13) $$\varsigma(-n,a) = -\frac{B_{n+1}(a)}{n+1}$$

With $a=1$ we obtain

(2.14) $$\varsigma(-n) = -\frac{B_{n+1}}{n+1}$$

and, in particular, we have for $n \ge 1$ the trivial zeros of the Riemann zeta function

(2.15) $$\varsigma(-2n) = 0$$

We may write (2.4) in the indeterminate form

(2.16) $$\varsigma(2n+1) = \frac{(2\pi)^{2n+1} B_{2n+1}}{2(2n+1)!\cos\left(\frac{[2n+1]\pi}{2}\right)}$$

Similarly, we may write (2.2) as

(2.17) $$\varsigma(s) = \frac{(2\pi)^s \beta(s)}{2s!\cos\left(\frac{\pi s}{2}\right)\cos[\pi(1-s)]}$$

and hence we have

$$\varsigma(2n+1) = \frac{(-1)^n (2\pi)^{2n+1}}{2(2n+1)!} \lim_{s\to 2n+1} \frac{\beta(s)}{\cos\left(\frac{\pi s}{2}\right)}$$

Using L'Hôpital's rule we have

(2.17.1) $$\varsigma(2n+1) = \frac{(-1)^n (2\pi)^{2n+1}}{2(2n+1)!} \lim_{s\to 2n+1} \frac{\beta'(s)}{\frac{\pi}{2}\sin\left(\frac{\pi s}{2}\right)} = \frac{2(-1)^n (2\pi)^{2n}}{(2n+1)!} \lim_{s\to 2n+1} \beta'(s)$$

Having regard to (2.8) and (2.11) let us consider

$$\beta(s) = -\cos[\pi(1-s)]B_s(1) = -\cos[\pi(1-s)]\sum_{m=0}^{\infty}\frac{1}{m+1}\sum_{k=0}^{m}\binom{m}{k}(-1)^k(1+k)^s$$

whereupon we have the derivative

$$\beta'(s) = -\cos[\pi(1-s)]\sum_{m=0}^{\infty}\frac{1}{m+1}\sum_{k=0}^{m}\binom{m}{k}(-1)^k(1+k)^s\log(1+k)$$

$$-\pi\sin[\pi(1-s)]\sum_{m=0}^{\infty}\frac{1}{m+1}\sum_{k=0}^{m}\binom{m}{k}(-1)^k(1+k)^s$$

We therefore have from (2.17.1)

(2.18) $$\varsigma(2n+1) = \frac{(-1)^n 2(2\pi)^{2n}}{(2n+1)!}\sum_{m=0}^{\infty}\frac{1}{m+1}\sum_{k=0}^{m}\binom{m}{k}(-1)^k(1+k)^{2n+1}\log(1+k)$$

Differentiating (2.9) we see that

$$(s-1)\varsigma'(s)+\varsigma(s) = -\sum_{m=0}^{\infty}\frac{1}{m+1}\sum_{k=0}^{m}\binom{m}{k}(-1)^k\frac{\log(1+k)}{(1+k)^{s-1}}$$

and therefore we get with $s=-2n$

$$(2n+1)\varsigma'(-2n) = \sum_{m=0}^{\infty}\frac{1}{m+1}\sum_{k=0}^{m}\binom{m}{k}(-1)^k(1+k)^{2n+1}\log(1+k)$$

since by (2.15) $\varsigma(-2n)=0$. Referring to (2.18) this then results in

(2.19) $$\varsigma(2n+1) = \frac{(-1)^n 2(2\pi)^{2n}}{(2n)!}\varsigma'(-2n)$$

which may also be obtained by differentiating Riemann's functional equation (2.5).

□

Referring to (2.2)

$$\beta(s) = 2\Gamma(s+1)\varsigma(s)(2\pi)^{-s}\cos\left(\frac{\pi s}{2}\right)\cos[\pi(1-s)]$$

we have

$$\log\beta(s) = \log 2+\log\Gamma(s+1)+\log\varsigma(s)-s\log(2\pi)+\log\cos\left(\frac{\pi s}{2}\right)+\log\cos[\pi(1-s)]$$

Logarithmic differentiation gives us

$$\frac{\beta'(s)}{\beta(s)} = \psi(s+1) + \frac{\varsigma'(s)}{\varsigma(s)} - \log(2\pi) - \frac{\pi}{2}\tan\left(\frac{\pi s}{2}\right) + \pi\tan[\pi(1-s)]$$

and we therefore have

$$\beta'(s) = \left[\psi(s+1) + \frac{\varsigma'(s)}{\varsigma(s)} - \log(2\pi) - \frac{\pi}{2}\tan\left(\frac{\pi s}{2}\right) + \pi\tan[\pi(1-s)]\right] 2\Gamma(s+1)\varsigma(s)(2\pi)^{-s}\cos\left(\frac{\pi s}{2}\right)\cos[\pi(1-s)]$$

With $s = 2n+1$ we obtain

(2.20) $$\beta'(2n+1) = (-1)^n \pi\,\Gamma(2n+2)\varsigma(2n+2)(2\pi)^{-2n-1}$$

which was also derived in a different manner in Woon's paper [13].

We recall the following expression for the generalised Stieltjes constants [3]

$$\gamma_p(u) = -\frac{1}{p+1}\sum_{m=0}^{\infty}\frac{1}{m+1}\sum_{k=0}^{m}\binom{m}{k}(-1)^k \log^{p+1}(u+k)$$

which gives us

$$\gamma_p = \gamma_p(1) = -\frac{1}{p+1}\sum_{m=0}^{\infty}\frac{1}{m+1}\sum_{k=0}^{m}\binom{m}{k}(-1)^k \log^{p+1}(1+k)$$

We previously showed in Section 5 of [3] how differentiating (2.10) may be used to determine a formula for the Stieltjes constants in terms of the higher derivatives of the Riemann zeta function $\varsigma^{(n)}(0)$.

□

Letting $s \to 1-s$ in (2.2) gives us

$$\beta(1-s) = 2\Gamma(2-s)\varsigma(1-s)(2\pi)^{s-1}\sin\left(\frac{\pi s}{2}\right)\cos(\pi s)$$

and using the Riemann functional equation (2.5) we get

$$\beta(1-s) = 2\Gamma(2-s)2(2\pi)^{-s}\Gamma(s)\cos\left(\frac{\pi s}{2}\right)\varsigma(s)(2\pi)^{s-1}\sin\left(\frac{\pi s}{2}\right)\cos(\pi s)$$

$$= \frac{1}{\pi}\Gamma(2-s)\Gamma(s)\varsigma(s)2\cos\left(\frac{\pi s}{2}\right)\sin\left(\frac{\pi s}{2}\right)\cos(\pi s)$$

$$= \frac{1}{\pi}\Gamma(2-s)\Gamma(s)\varsigma(s)\sin(\pi s)\cos(\pi s)$$

Since $\Gamma(2-s) = \Gamma(1+1-s) = (1-s)\Gamma(1-s)$ we have

$$= \frac{1}{\pi}(1-s)\Gamma(s)\Gamma(1-s)\varsigma(s)\sin(\pi s)\cos(\pi s)$$

and therefore using Euler's reflection formula we obtain

(2.21) $\beta(1-s) = (1-s)\varsigma(s)\cos(\pi s)$

With $s = -1$ we get

$$\beta(2) = -2\varsigma(-1) = \frac{1}{6}$$

With $s = 1/2$ we get

$$\beta(1/2) = 0$$

It seems possible that Carlson's theorem [12, p.186] may be applicable to (2.2). Carlson's theorem applies to functions $f(z)$ analytic for $\mathrm{Re}(z) \geq 0$ satisfying the bound $|f(z)| = O\left(e^{\mu|z|}\right)$, $\mu < \pi$. It asserts that if $f(z) = 0$ for $z = 0,1,2...$, then $f(z) = 0$ identically.

Donal F. Connon
Elmhurst
Dundle Road
Matfield
Kent TN12 7HD
dconnon@btopenworld.com